\input ./style/arxiv-general.cfg
\documentclass[aap,MSNbibl,dvips]{arximspdf}
\makeatletter
   \@ifpackageloaded{graphicx}{}{\usepackage{graphicx}}
\makeatother

\doi{10.1214/15-AAP1131}
\volume{26}
\issue{3}
\pubyear{2016}
\firstpage{1727}
\lastpage{1742}
\docsubty{FLA}

\makeatletter
\newcommand{\lleft}{\left}
\newcommand{\rright}{\right}
\newcommand{\eqref}[1]{(\ref{#1})}
\newcommand{\RR}{\mathbb{R}}
\newcommand{\GG}{\mathbb{G}}
\renewcommand{\SS}{\mathbb{S}}
\newcommand{\PP}{\mathbb{P}}
\newcommand{\EE}{\mathbb{E}}
\newcommand{\vd}{ \mathrm{d}}

\newcommand{\dd}{\mathrm{d}}
\newcommand{\cN}{\mathcal{N}}
\newcommand{\cU}{\mathcal{U}}
\newcommand{\cC}{\mathcal{C}}
\newcommand{\cCc}{\mathcal{C}_{\mathrm{c}}}
\newcommand{\cB}{\mathcal{B}}

\newcommand{\cF}{\mathcal{F}}

\newcommand{\ft}{\mathfrak{t}}
\newcommand{\fl}{\mathfrak{l}}
\newcommand{\fr}{\mathfrak{r}}
\newcommand{\fs}{\mathfrak{s}}
\newcommand{\ff}{\mathfrak{f}}
\newcommand{\fh}{\mathfrak{h}}
\newcommand{\fk}{\mathfrak{k}}
\newcommand{\fg}{\mathfrak{g}}
\newcommand{\fz}{\mathfrak{z}}
\newcommand{\fu}{\mathfrak{u}}
\newcommand{\cZ}{\mathcal{Z}}
\newcommand{\cL}{\mathcal{L}}
\newcommand{\cA}{\mathcal{L}}
\newcommand{\Dom}{\operatorname{Dom}}
\newcommand{\grandO}{\mathrm{O}}
\newcommand{\cIG}{\mathcal{IG}}
\newcommand{\uL}{\mathrm{L}}
\newcommand{\uH}{\mathrm{H}}

\newcommand{\equaldist}{\stackrel{\mathrm{dist}}{=}}
\newcommand{\sgn}{\operatorname{sgn}}

\newcommand{\Gel}{G^{\mathrm{e}}}
\newcommand{\gel}{g^{\mathrm{e}}}
\newcommand{\Pel}{P^{\mathrm{e}}}
\newcommand{\Pre}{P^{\mathrm{r}}}
\newcommand{\Gre}{G^{\mathrm{r}}}

\newtheorem{theorem}{Theorem}
\newtheorem{proposition}{Proposition}

\newproclaim{definition}{Definition}
\newproclaim{remark}{Remark}
\makeatother

\begin{document}
\begin{frontmatter}

\title{The snapping out Brownian motion}
\runtitle{The snapping out Brownian motion}

\begin{aug}
\author[A]{\fnms{Antoine} \snm{Lejay}\corref{}\thanksref{T1}\ead[label=e1]{Antoine.Lejay@univ-lorraine.fr}}
\thankstext{T1}{Supported by the
ANR SIMUDMRI (ANR-10-COSI-SIMUDMRI).}

\runauthor{A. Lejay}
\affiliation{Inria Nancy Grand-Est}
\address[A]{Institut \'{E}lie Cartan de Lorraine\\
UMR 7502\\
Universit\'{e} de Lorraine\\
Vand\oe uvre-l\`{e}s-Nancy, F-54500\\
France\\
and\\
Institut \'{E}lie Cartan de Lorraine\\
UMR 7502\\
CNRS\\
Vand\oe uvre-l\`{e}s-Nancy, F-54500\\
France\\
and\\
TOSCA\\
Inria\\
Villers-l\`{e}s-Nancy, F-54600\\
France\\
\printead{e1}}

\end{aug}

%
\received{\smonth{1} \syear{2013}}
%
\revised{\smonth{7} \syear{2015}}

%
\begin{abstract}
We give a probabilistic representation of a one-dimensional diffusion
equation where the solution is discontinuous at $0$
with a jump proportional to its flux. This kind of interface
condition is usually seen as a semi-permeable barrier.
For this, we use a process called here the snapping out Brownian motion,
whose properties are studied. As this construction is motivated
by applications, for example, in brain imaging or in chemistry,
a simulation scheme is also provided.
\end{abstract}

%
\end{frontmatter}

\section{Introduction}

Many diffusion phenomena have to deal with interface conditions.
Let $D$ be a diffusivity coefficient which is smooth away
from a regular surface $S$, but presents some discontinuity
there. In this case, the solution to the diffusion equation
%
%
\begin{equation}
\label{eq-pde1} \partial_t u(t,x)=\tfrac{1}{2}\nabla\bigl(D(x)
\nabla u(t,x)\bigr)=0 \qquad\mbox{with }u(0,x)=f(x)
\end{equation}
has to be understood as a \emph{weak solution}. However,
$u$ is smooth away from $S$ and satisfies
%
%
\begin{eqnarray}
\label{eq-pde2} u(t,x+)&=&u(t,x-)\quad\mbox{and}
\nonumber
\\[-8pt]
\\[-8pt]
\nonumber
 D(x+)n^+(x)\cdot\nabla u(t,x+)&=&
D(x-)n^-(x)\cdot\nabla u(t,x-),
\end{eqnarray}
for $x\in S$, when $S$ is assumed to separate
locally $\RR^d$ into a ``$+$'' and a ``$-$'' part and
where $n^\pm$ is a vector normal to $S$ at $x$ pointing
to the ``$\pm$'' side.
The second condition is called the \emph{continuity of the flux}.

Now, let us assume that $D$
takes scalar values, and is constant away from a thin layer
of width $2\ell$ enclosed between two parallel surfaces $S^+$ and $S^-$.
When the width $\ell$ of the layer tends to $0$, $S^+$ and $S^-$ merge
into a single interface located on a surface $S$.

\begin{figure}

\includegraphics{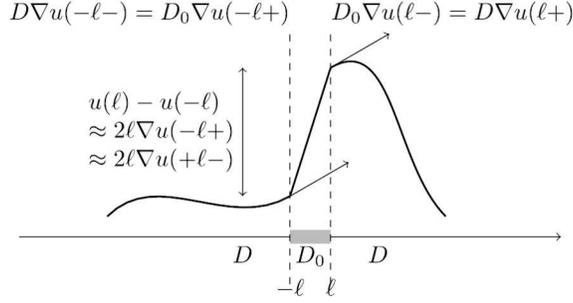}

\caption{The thin layer problem.}
\label{fig-1}
\end{figure}

When the diffusivity $D_0$ decreases to $0$ with $\ell$
and $D_0/\ell\to\lambda>0$,
then the solution to \eqref{eq-pde1} converges to a function $v$
satisfying \eqref{eq-pde1} away from $S$ with
the interface condition for $x\in S$:
%
%
\begin{eqnarray}
\label{eq-pde3} \nabla v(t,x+)&=&\nabla v(t,x-)\quad\mbox{and}
\nonumber
\\[-8pt]
\\[-8pt]
\nonumber
 \frac{\lambda}{2}
\bigl(v(t,x+)-v(t,x-)\bigr)&=&D(x\pm)\nabla v(t,x\pm).
\end{eqnarray}
The solution has a continuous flux on $S$
but is discontinuous on $S$ (see, e.g.,
\cite{sanchez}, Chapter~13). A heuristic explanation is
given Figure~\ref{fig-1}.

%

If $D$ is smooth on $\RR^d$, it is well known that
%
%
\begin{equation}
\label{eq-intro-2} u(t,x)=\EE_x\bigl[f(X_t)\bigr],
\end{equation}
where $X$ is the diffusion process generated by $\frac{1}{2}\nabla
(D\nabla)$
which is solution under $\PP_x$
to the stochastic differential equation (SDE)
%
%
\begin{equation}
\label{eq-intro-1} X_t=x+\int_0^t
\sigma(X_s)\,\vd B_t+\int_0^t
\frac{1}{2}\sum_{i=1}^d
\frac{D_{i,\cdot}}{\partial x_i}(X_s)\,\vd s\qquad \mbox{with }\sigma\sigma
^{\mathrm{T}}=D
\end{equation}
for a Brownian motion $B$.

When $D$ presents some discontinuities,
\eqref{eq-intro-1} has no longer a meaning.
However, a Feller processes $(X,(\cF_t)_{t\geq0},(\PP_x)_{x\in\RR})$ is
associated to $\frac{1}{2}\nabla(D\nabla\cdot)$
for which \eqref{eq-intro-2} holds. In particular,
the marginal distributions $X_t$ have a density
$p(t,x,\cdot)$ under~$\PP_x$, where $p(t,x,y)$
is the fundamental solution to \eqref{eq-pde1} (see, e.g., \cite{stroock}).

Let us now assume that the dimension of the space
is equal to $1$ and that $D$ is discontinuous
at some separated points $\{x_i\}$ with
left and right limit there, and smooth elsewhere.
The process $X$ is
solution to a SDE with local time. The It\^{o}--Tanaka
formula is the key tool to manipulate it, and
several simulation algorithms have been proposed
(see the references in \cite{lejay-pichot}, e.g.).
The process called the \emph{Skew Brownian motion} is the
main tool for this construction \cite{lejay-sbm,lejay-martinez}.

Coming back to the thin layer problem, we assume that $D$
is constant and equal to $D_1$ on $(-\infty,-\ell)$
and $(\ell,\infty)$, and to $D_0$ on $(-\ell,\ell)$.
The associated stochastic process is
solution to
\[
X_t=x+\int_0^t
\sqrt{D(X_s)}\,\vd B_s+\frac{D_1-D_0}{D_1+D_0}L_t^{\ell}(X)
+\frac{D_0-D_1}{D_1+D_0}L_t^{-\ell}(X),
\]
where $L^{\pm\ell}_t(X)$ is the local time of $X$ at $\pm\ell$
\cite{lejay-martinez}.

Letting $D_0/\ell$ converging to $2\kappa$ with $\ell\to0$,
one may expect that $X$ converges in distribution to a
stochastic process $Y$ such that the solution
to \eqref{eq-pde1} with the interface condition \eqref{eq-pde3}
is given by $v(t,x)=\EE_x[f(Y_t)]$.

The article then aims at constructing and giving several properties
related to the process $Y$ which we call
a \emph{snapping out Brownian motion} (SNOB).
This process is Feller on $\GG=(-\infty,0-]\cup[0+,+\infty)$
but not on $\RR$. The intervals
in the definition of $\GG$ are disjoint
so that $0$ corresponds either to $0+$ or $0-$
seen as distinct points.

The behavior of this process is the following:
Assume that its starting point is $x\geq0$.
It behaves as a positively reflected Brownian motion
until its local time is greater than an independent exponential random variable
of parameter $2\kappa$. Then its decides its sign
with probability $1/2$ and starts afresh as a new reflected Brownian motion,
until its local time is greater than a new exponential random variable,
and so on. Using the properties of the exponential
random variable, it is equivalent to assert that the particle
changes its sign when its local time is greater
than an exponential random variable with parameter $\kappa$,
and behaves like a positively or negatively reflected
Brownian motion between these switching times.

Its name is justified by the following fact:
As the time at which the particle possibly changes
it signs is the same as for the \emph{elastic
Brownian motion} \cite{feller,ito-mckean,karlin,grebenkov}
(also called the \emph{partially reflected Brownian motion}),
it could also be seen as some elastic Brownian motion
which is reborn once killed.

The elastic Brownian motion, also
called a \emph{partially reflected Brownian motion},
is associated to the Robin boundary condition and
has then many applications \cite{grebenkov,singer,erban}.
This process is the ``basic brick'' for constructing
the SNOB.

The behavior of the SNOB
justifies also the old heuristic that the interface
condition \eqref{eq-pde3} corresponds to a
\emph{semi-permeable barrier}, which arises, for example,
in diffusion magnetic resonance imaging \cite{fiermans}
or in chemistry \cite{andrews,erban}. The interface condition
\eqref{eq-pde3} is different from \eqref{eq-pde2}, to which
is associated a Skew Brownian motion and where the particle
crosses the interface when it reaches it, and which corresponds
to a \emph{permeable barrier} (see references in \cite{lejay-sbm,lejay-pichot}).

Here, we work under the condition of a single interface at $0$.
In short time, it is sufficient to describe the
behavior of the process even in a more complex media, since other
interface or boundary conditions far enough have
``exponentially small'' influence on
the distribution of the process. This is sufficient
for simulation purposes, where particles positions
are represented by the stochastic process and move
according to its dynamic.

Using similar computations, one may generalize our
work to the case where $D(x)=D^+$ if $x\geq0$,
$D^-$ if $x\leq0$ and an interface condition
\[
\nabla u(t,0+)=\beta\nabla u(t,0-)\quad\mbox{and}\quad \lambda u(t,0+)-\mu u(t,x-)=
\nabla u(t,x+)
\]
with $\lambda,\mu>0$. Diffusions on graphs
specified by a condition at each vertex could
also be considered, which could be of interest
in several applications. This process has been described
without proof by Bobrowski in \cite{bobrowski}, which
have studied its limit behavior when the diffusion coefficients
increase.

Although the SNOB may be seen
as a diffusion on a graph, it is not a diffusion on a metric graph,
where the edges are joined by vertices.
Such diffusions have been classified by Freidlin and Wentzell
in \cite{freidlin2,freidlin}. The conditions that are
required at the vertices of the graphs are some
extension of the possible boundary conditions
for a Markov process studied by Feller \cite{feller}.
See also \cite{kopytko}, for example, for the related problem of
pasting diffusions.%
\setcounter{footnote}{1}\footnote{The article \cite{portenko} defines a notion
of \emph{semipermeable membrane} which is different from ours,
where the solution is continuous with a discontinuous gradient.}

Our interface condition does not fall in these categories.
Our process is best thought as a kind of \emph{random evolution
process} which switches
back and forth randomly among a collection of processes
(see, e.g., \cite{griego,siegrist}).

\textit{Outline}.
In Sections \ref{sec:ebm} and \ref{sec:piecing}, we present quickly the main
results related to the elastic Brownian motions and the piecing out procedure.
The SNOB is constructed in Section~\ref{sec:snob}
through its resolvent. In Section~\ref{sec:layer}, we show the relationship
between the SNOB and the thin layer problem. Finally, in Section~\ref
{sec:simulation},
we show how to simulate this process.

\section{Elastic Brownian motion}
\label{sec:ebm}

Let $(R_t)_{t\geq0}$ a reflected Brownian motion, and denote by
$(L_t)_{t\geq0}$ its
symmetric local time at $0$. We add a cemetery point $\dag$ to $\RR_+$.
For a constant $\kappa>0$, we consider an exponential random
variable $\xi$ with parameter $\kappa$ independent from $B$. Set
\[
Z_t=\cases{ R_t,&\quad$\mbox{if }L_t\leq\xi$,
\vspace*{2pt}
\cr
\dag,&\quad$\mbox{if }L_t>\xi$.}
\]
Thanks to the properties of the local time, this process,
called the \emph{elastic Brownian motion} (EBM),
is still a strong Markov process. Its semi-group is
\[
\Pel_t f(x)=\EE_x\bigl[\exp(-\kappa L_t)f(X_t)
\bigr]
\]
for $f$ in the set $\cC_0(\RR_+,\RR)$ of continuous
functions that vanishes at infinity. Closed form
expressions of the density transition function are given in \cite
{singer,gallavotti}.

Let $\fk$ be the time at which the EBM is killed, which
means $\fk=\inf\{t>0 | L_t\geq\xi\}$. This is a stopping time.
Since
the local time increases only on the closure of $\cZ=\{t>0 | X_t=0\}$,
it holds that $Z_{\fk}=0$ almost surely.
Using standard computations in the inverse of the local time
of the Brownian motion,
%
%
\begin{equation}
\label{eq:ebm:2} \psi(x,\alpha)=\EE_x\bigl[\exp(-\alpha\fk)\bigr]=
\frac{\kappa}{\sqrt{2\alpha}+\kappa}\exp(-\sqrt{2\alpha}x).
\end{equation}

Using the It\^{o} formula, it is easily shown
that $u(t,x)=\Pel_t f(x)$ is solution to the heat
equation with Robin (or third kind) boundary condition \cite
{grebenkov,papanicolaou,bass08a}
\[
\cases{ \displaystyle\frac{\partial u(t,x)}{\partial t} =\frac{1}{2}\triangle
u(t,x),&\quad$\mbox{on }(0,+
\infty)^2,$\vspace*{2pt}
\cr
\displaystyle\frac{\partial u(t,0)}{\partial x}=\kappa u(t,0). }
\]

For a Markov process $X$, let us recall that its resolvent
$(G_\alpha)_{\alpha>0}$ is a family of operators
defined by
$G_\alpha f(x)=\EE_x [\int_0^{+\infty} e^{-\alpha s}f(X_s)\,\vd s ]$
for any $f\in\cC_0$ and any $\alpha>0$.
It has a density $g_\alpha$ when $G_\alpha f(x)=\int g_\alpha
(x,y)f(y)\,\vd y$.

Using standard computations on the Green functions,
the density $\gel_\alpha(x,y)$ of the resolvent of the EBM
is for $x,y\geq0$,
\[
\gel_\alpha(x,y)= \frac{1}{\sqrt{2\alpha}} \cases{ \displaystyle\frac{\sqrt{2\alpha
}-\kappa}{\sqrt{2\alpha}+\kappa}
e^{-\sqrt{2\alpha}(y+x)} +e^{-\sqrt{2\alpha}(x-y)}, &\quad$\mbox{for }y\in
[0,x]$,\vspace*{2pt}
\cr
\displaystyle e^{\sqrt{2\alpha}(x-y)}+\frac{\sqrt{2\alpha}-\kappa}{\sqrt{2\alpha
}+\kappa}e^{-\sqrt{2\alpha}(x+y)}, &\quad $\mbox{for }y\geq x$.}
\]

We extend the EBM to a process on $\GG$ by symmetry, so
that its resolvent becomes
%
%
\begin{equation}
\label{eq:ebm:1} \Gel_\alpha f(x):=\EE_x \biggl[\int
_0^\fk e^{-\alpha s}f(X_s)\,\vd s
\biggr] =\int_0^{+\infty} \gel_\alpha\bigl(|x|,y\bigr)f
\bigl(\sgn(x)y\bigr)\,\vd y
\end{equation}
for $x\in\GG$. This process evolves either on $\RR_-$ or $\RR_+$
but never crosses $0$ and is naturally identified with a process on $\GG$.

\section{Piecing out Markov processes}
\label{sec:piecing}

The procedure of piecing out is a way to construct
a Markov process from a killed one. We present in this section a result
due to Ikeda,
Nagasawa and Watanabe \cite{ikeda} (similar considerations
are given in~\cite{meyer}).

On a probability space $(\Omega,\cF,\PP)$ and a state space $\SS$,
let $((X_t)_{t\geq0},(\cF_t)_{t\geq0},\break (\PP_x)_{x\in\overline{\SS}})$
be a right continuous strong Markov process living
in the extended state space $\SS^\dag=\SS\cup\{\dag\}$ with a death
point $\dag$.
The lifetime of $X$ is denoted by $\fk$.

The shift operator associated
to $X$ is denoted by $(\theta_t)_{t\geq0}$.

We also consider a family $\mu$ defined
on $\Omega\times\SS^\dag$ such that $\mu(\omega,\cdot)$
is a probability measure on $\SS^\dag$ and for any fixed Borel
subset $A$, $\mu(\cdot,A)$ is $\sigma(X_t,t\geq0)$-measurable.
We assume additionally that
$\mu(\omega,\vd y)=\delta_\dag(\dd y)$ when
$\fk(\omega)=0$ and
\[
\PP_x\bigl[\mu(\omega,\vd y)=\mu(\theta_{\ft(\omega)}\omega,\vd y),
\ft(\omega)<\fk(\omega)\bigr] =\PP_x[\ft<\fk]
\]
for any stopping time $\ft$.
The family $\mu$, called an \emph{instantaneous distribution},
describes the way the process
is reborn once killed.

Let $\widehat{\Omega}$ be the product of an infinite,
countable, number of copies of $\Omega\times\SS^\dag$.
We define $X$ on $\widehat{\Omega}$ by
\[
X_t(\widehat{\omega}) =\cases{ x_t(\omega_1),&\quad $
\mbox{if }t\in\bigl[0,\fk(\omega_1)\bigr)$,\vspace*{2pt}
\cr
y_1, &\quad $\mbox{if }t=\fk(\omega_1)$,\vspace*{2pt}
\cr
x_{t-\fk(\omega_1)}(
\widetilde{\omega}_2),&\quad $\mbox{if } t\in\bigl(\fk(\omega_1),
\fk(\omega_1)+\fk(\omega_2)\bigr)$,\vspace*{2pt}
\cr
y_2,&$\quad \mbox{if }t=\fk(\omega_2)$,\vspace*{2pt}
\cr
\cdots&\vspace*{2pt}
\cr
\dag,&\quad $\mbox{ if }t\geq\fk(\omega_1)+\cdots+
\fk_N(\omega_N)$ }
\]
with $\widehat{\omega}=(\omega_1,y_1,\omega_2,y_2,\ldots)\in
\widehat{\Omega}$ and $N=\inf\{k\geq0;\fk(\omega_k)=0\}$.

We consider the probability measure
\begin{eqnarray*}
&&\widehat{\PP}_x\bigl[\vd\omega_1,\vd x^1,
\ldots,\vd\omega_n,\vd x^n\bigr]
\\
&&\qquad=\PP_x\bigl[\vd\omega^1\bigr]\mu\bigl(
\omega^1,\vd x^1\bigr) \PP_{x^1}\bigl[\vd
\omega^2\bigr]\mu\bigl(\omega^1,\vd x^2\bigr)
\cdots\PP_{x^n}\bigl[\vd\omega^2\bigr]\mu\bigl(
\omega^n,\vd x^n\bigr).
\end{eqnarray*}
Under this measure $\widehat{\PP}_x$,
when the path $X(\omega)$ is killed, we let it reborn
by placing it at the point $x_1$ with probability $\mu(\omega,\vd x_1)$
and then start again.

We left the technical details about the construction of the probability
space and the filtration and presents the main result on piecing out
Markov process.

%
\begin{theorem}[(\cite{ikeda})]
\label{thm-piecing}
Using the above defined notation, there exists a probability
space $(\widehat{\Omega},\widehat{\cB},\widehat{\PP})$ and a filtration
$(\widehat{\cB}_t)_{t\geq0}$ on which
$(X,(\widehat{\cB}_t)_{t\geq0},(\widehat{\PP}_x)_{x\in\SS^\dag})$
is a strong Markov process on $\SS^\dag$ with
$\PP_\dag[X_t=\dag, \forall t\geq0]=1$.
\end{theorem}

\section{The snapping out Brownian motion}
\label{sec:snob}

%
\begin{definition} A \emph{snapping out Brownian motion} (SNOB) $X$
is a strong Markov stochastic process living on $\GG$
constructed
by making EBM reborn on $0+$ or $0-$ with probability $1/2$
using the piecing-out procedure.
\end{definition}

The sign of $X$ changes with probability $1/2$ when
its local time $L_t$ at $0$ is greater than
$\fu_k$ with $\fu_0=0$, $\fu_k-\fu_{k-1}\sim\exp(\kappa)$
is independent from $(\fu_i)_{i\leq k-1}$. From the properties
of the exponential and binomial distributions, the sign of $X$ changes
when its local time is greater than
$\fs_k$ with $\fs_0=0$, $\fs_k-\fs_{k-1}\sim\exp(\kappa/2)$
is independent from $(\fs_i)_{i\leq k-1}$.

It is also immediate that $|X|$ is a reflected Brownian motion,
where $|\cdot|$ is the canonical projection of $\GG$ onto $[0,+\infty)$.

%
\begin{proposition}
\label{prop:main}
The resolvent family $(G_\alpha)_{\alpha>0}$ of the SNOB
is solution to
\begin{eqnarray}
\biggl(\alpha-\frac{1}{2}\triangle\biggr)G_\alpha f(x)=f(x)\qquad\mbox{
for }x\in\GG
\nonumber\\
\eqntext{\mbox{with } \nabla G_\alpha f(0+)=\nabla G_\alpha f(0-)\mbox{
and } \displaystyle\frac{\kappa}{2}\bigl(G_\alpha f(0+)-G_\alpha f(0-)
\bigr)=\nabla G_\alpha f(0)}
\end{eqnarray}
for any bounded, continuous function $f$ on $\GG$ that vanishes at infinity.
\end{proposition}

This proposition identifies the infinitesimal generator
of the process $X$. The points $0+$ and $0-$ are then interpreted
as the sides of a semi-permeable barrier.

\begin{pf*}{Proof of Proposition~\protect\ref{prop:main}}
From this very construction and the strong Markov property,
for any continuous function $f$ on $\GG$ which vanishes at infinity,
%
%
\begin{equation}
\label{eq:1} G_\alpha f(x) =\Gel_\alpha f(x) +
\frac{\psi(|x|,\alpha)}{2} \bigl(G_\alpha f(0+)+G_\alpha f(0-) \bigr),
\end{equation}
where $\Gel_\alpha$ is defined by \eqref{eq:ebm:1}.

Using $x=0+$ and $x=0-$ in \eqref{eq:1} and summing
the two resulting equations leads to
%
%
\begin{eqnarray}
\label{eq:1bis} G_\alpha f(x)=\Gel_\alpha f(x)+\frac{\kappa e^{-\sqrt
{2\alpha}|x|}}{2\sqrt{2\alpha}}
\beta(f)
\nonumber
\\[-8pt]
\\[-8pt]
\eqntext{\mbox{with }\beta(f)=\Gel_\alpha f(0+)+\Gel_\alpha
f(0-).}
\end{eqnarray}
Then
%
%
\begin{eqnarray}
\label{eq:3} G_\alpha f(x)+G_\alpha f(-x)&=&\Gel_\alpha
f(x)+\Gel_\alpha f(-x)+ \frac{\kappa}{\sqrt{2\alpha}}e^{-\sqrt{2\alpha
}|x|}\beta(f),
\\
\label{eq:4} G_\alpha f(x)-G_\alpha f(-x)&=&\Gel_\alpha
f(x)-\Gel_\alpha f(-x).
\end{eqnarray}
Derivating \eqref{eq:3} and setting $x=0+$,
since $\nabla\Gel_\alpha f(0\pm)=\pm\kappa\Gel_\alpha f(0\pm)$,
\[
\nabla G_\alpha f(0+)-\nabla G_\alpha f(0-)=0.
\]
Derivating \eqref{eq:4},
\begin{eqnarray*}
2 \nabla G_\alpha f(0\pm)& =&\nabla G_\alpha f(0+)+\nabla
G_\alpha f(0-) =\nabla\Gel_\alpha f(0+)+\nabla
\Gel_\alpha f(0-)
\\
&=&\kappa\bigl(\Gel_\alpha f(0+)-\Gel_\alpha f(0-)\bigr) =\kappa
\bigl(G_\alpha f(0+)-G_\alpha f(0-)\bigr).
\end{eqnarray*}

In addition, it is easily seen that
$ (\alpha-\frac{1}{2}\triangle)G_\alpha f=f$
since $\psi(x,\alpha)$
is solution to $ (\alpha-\frac{1}{2}\triangle) \psi(x,\alpha)=0$.
The resolvent is then identified.
\end{pf*}

%
\begin{proposition}
\label{prop:representation}
The semi-group $(P_t)_{t\geq0}$ of
the SNOB has the following representation:
%
%
\begin{eqnarray}
\label{eq:5} P_t f(x)&=& \EE_x \biggl[ \biggl(
\frac{1+e^{-\kappa L_t}}{2} \biggr) f\bigl(\sgn(x)|B_t|\bigr) \biggr]
\nonumber
\\[-8pt]
\\[-8pt]
\nonumber
 &&{}+
\EE_x \biggl[ \biggl(\frac{1-e^{-\kappa L_t}}{2} \biggr) f\bigl(-
\sgn(x)|B_t|\bigr) \biggr]
\end{eqnarray}
for a Brownian motion $B$.
\end{proposition}

\begin{pf}
Let us decompose a function $f$ as its even and odd parts:
\[
\hat{f}(x)=\tfrac{1}{2}\bigl(f(x)+f(-x)\bigr)\quad \mbox{and}\quad \check{f}(x)=
\tfrac{1}{2}\bigl(f(x)-f(-x)\bigr).
\]
Then $\Gel_\alpha\hat{f}(-x)=\Gel_\alpha\hat{f}(x)$
and $\Gel_\alpha\check{f}(-x)=-\Gel_\alpha\check{f}(x)$,
so that $\beta(\check{f})=0$ for $\beta$ defined by \eqref{eq:1bis}.
Thus $G_\alpha\check{f}(x)=\Gel_\alpha\check{f}(x)$.
In addition,
since $\hat{f}(|x|)=\hat{f}(x)$ and the SNOB has the same
distribution as the reflected Brownian motion $|B|$,
\[
G_\alpha\hat{f}(x)=\Gre_\alpha\hat{f}(x):=\EE_{x}
\biggl[ \int_0^{+\infty} e^{-\alpha s}
\hat{f}\bigl(|B_s|\bigr)\,\vd s \biggr].
\]
This gives an alternative representation for the resolvent
of the SNOB:
$G_\alpha f(x)=\Gre_\alpha\hat{f}(x)
+\Gel_\alpha\check{f}(x)$.
Inverting the resolvent to recover the semi-group
$(P_t)_{t\geq0}$,
\[
P_t f(x)=\Pre_t \hat{f}(x)+\Pel_t
\check{f}(x) =\EE_x\bigl[\hat{f}\bigl(|B_t|\bigr)\bigr]+
\EE_x\bigl[\exp(-\kappa L_t)\check{f}\bigl(
\sgn(x)|B_t|\bigr)\bigr].
\]
This expression could be arranged as \eqref{eq:5}.
\end{pf}

\section{The thin layer problem}

\label{sec:layer}

We now fix $\varepsilon>0$ and we consider the process $X^\varepsilon$
generated by (see, e.g., \cite{stroock} for general considerations
on this process)
\[
\cL^\varepsilon:=\frac{1}{2}\frac{\partial}{\partial x} \biggl(
a^\varepsilon(x) \frac{\partial}{\partial x} \biggr)\qquad \mbox{with
}a^\varepsilon(x):=
\cases{1,&\quad$\mbox{when }x\notin[-\varepsilon,\varepsilon]$,\vspace*{2pt}
\cr
\kappa
\varepsilon,&\quad$\mbox{when }x\in[-\varepsilon,\varepsilon]$}
\]
whose domain $\Dom(\cL^\varepsilon)=\{f\in\uL^2(\RR) | \cL^\varepsilon
f\in\uL^2(\RR)\}$
is a subset of the Sobolev space $\uH^1(\RR)$ [hence, any function in
$\Dom(\cL^\varepsilon)$
is identified with a continuous function], where
$\uL^2(\RR)$ is the set of square integrable functions on $\RR$
with scalar product $\langle f,g\rangle=\int_{\RR} f(x)g(x)\,\vd x$.
Let us set $[h](x):=h(x-)-h(x+)$ and
%
%
\begin{equation}
\label{eq:Deps} D^\varepsilon:=\lleft\{ f\in\cC^2\bigl((-
\infty,-\varepsilon)\cup(-\varepsilon,\varepsilon)\cup(\varepsilon
,\infty)\bigr)
\middle|
\begin{array} {l} f,f''\in
\uL^2(\RR),
\\
{}[f](\pm\varepsilon)=0,
\\
\bigl[a^\varepsilon\nabla f\bigr](\pm\varepsilon)=0
\end{array}
\rright\}.
\end{equation}
For $k\geq0$, we write $\cCc^k(\RR)$ the set of functions with compact support
and continuous derivatives up to order $k$.
With an integration by parts, for $f\in D^\varepsilon$ and $g\in\cCc
^2(\RR)$,
\begin{eqnarray*}
\bigl\langle(\alpha-L)f,g\bigr\rangle&= &\alpha\langle f,g\rangle+\int
_{\RR} a^\varepsilon(x)\nabla f(x)\nabla g(x)\,\vd x\\
&&{} +
\bigl[a^\varepsilon\nabla f\bigr](-\varepsilon)g(-\varepsilon)-
\bigl[a^\varepsilon\nabla f\bigr](\varepsilon)g(\varepsilon).
\end{eqnarray*}
Using this formula and the regularity of the solution to $(\alpha-L)f=g$
when $g\in\cC^\infty(I,\RR)$ with $-\varepsilon,\varepsilon\notin I$,
we easily get that $D^\varepsilon$ contains\vspace*{1pt} $(\alpha-\cL^\varepsilon
)^{-1}(\cCc^\infty(\RR))$
and is then dense in $\Dom(\cL^\varepsilon)$ for the operator norm
$(\langle f,f\rangle+\langle Lf,Lf\rangle)^{1/2}$.

A fundamental solution may be associated to $\cL^\varepsilon$,
as well as a resolvent density $g_\alpha^\varepsilon$, which we will
compute explicitly.

This operator is self-adjoint with respect to $\langle\cdot,\cdot
\rangle$,
so that its resolvent density satisfies $g_\alpha^\varepsilon
(x,y)=g_\alpha^\varepsilon(y,x)$.
This process is a Feller process, and is a strong solution
to the SDE with local time
\[
X^\varepsilon_t=x+\int_0^t
\sqrt{a^\varepsilon\bigl(X^\varepsilon_s\bigr)}\,\vd
B_s +\eta_\varepsilon L_t^{\varepsilon}
\bigl(X^\varepsilon\bigr) -\eta_\varepsilon L_t^{-\varepsilon}
\bigl(X^\varepsilon\bigr) \qquad\mbox{with } \eta_\varepsilon=\frac{1-\kappa
\varepsilon}{1+\kappa\varepsilon},
\]
where $B$ is a Brownian motion and $L_t^x(X^\varepsilon)$
is the symmetric local time at $x$ of $X^\varepsilon$
(see, e.g., \cite{lejay-martinez}, and
\cite{legall,bass} among others for general results on SDEs with local time).

In \cite{feller}, Section~11, the elastic Brownian motion is
constructed as the limit of a process
which either jumps at $\varepsilon$
or is killed with probability $\kappa\varepsilon$ when
it arrives at~$0$.

Using the piecing out procedure, we construct a strong Markov process
$Z^\varepsilon$
by considering the process $X^\varepsilon$ which is instantaneously
replaced at $-\varepsilon$ or $\varepsilon$ with probability $1/2$
when it reaches $0$, and then behaving again as $X^\varepsilon$ until
it reaches $0$, and so on.
This process $Z^\varepsilon$ could be identified as a process living in
$\GG$ by defining
$\PP_{0+}$ as $\PP_{\varepsilon}$ and $\PP_{0-}$ as $\PP_{-\varepsilon}$,
since the process is instantaneously killed when at $0$.

%
\begin{theorem}
\label{thm:main}
The process $Z^\varepsilon$ with $Z^\varepsilon_0=x$ converges in distribution
to the SNOB starting from $x$ in the Skorohod topology.
\end{theorem}

The proof relies on the next two results.

%
\begin{proposition}\label{prop3}
Let $g^\varepsilon_\alpha$ be the resolvent density of $X^\varepsilon$.
Then $g^\varepsilon_\alpha(x,y)$ converges to $g(x,y)$
for any $x,y\neq0$ and any $\alpha>0$ as $\varepsilon\to0$.
\end{proposition}

%
\begin{remark}
This result follows from classical results in deterministic homogenization
theory (see, e.g., \cite{sanchez}) where the convergence
holds in Sobolev spaces. Here, we consider a direct computational
proof for the convergence of the Green kernel, which we use later.
\end{remark}

\begin{pf*}{Proof of Proposition \ref{prop3}}
We assume that $x>0$ and we set $\mu:=\sqrt{2\alpha}$
for some $\alpha>0$. The resolvent density $g^\varepsilon_\alpha$ of
$X^\varepsilon$ has the form,
for $x>\varepsilon$,
\[
g^\varepsilon_\alpha(x,y)=\cases{ C_\varepsilon(x)e^{-\mu y},&\quad $
\mbox{for $y>x$},$\vspace*{2pt}
\cr
A_\varepsilon(x)e^{-\mu y}+B_\varepsilon(x)e^{\mu y},&\quad $
\mbox{for }y\in[\varepsilon,x],$\vspace*{2pt}
\cr
H_\varepsilon(x)e^{\mu y/\sqrt{\kappa\varepsilon}}+E_\varepsilon
(x)e^{-\mu y/\sqrt{\kappa\varepsilon}},&\quad $
\mbox{for }y\in[-\varepsilon,\varepsilon],$\vspace*{2pt}
\cr
F_\varepsilon(x)e^{\mu y},&\quad $
\mbox{for }y<-\varepsilon.$}
\]
By this, we mean that for any bounded, measurable function $f$,
\[
\EE_x \biggl[\int_0^{+\infty}
e^{-\alpha t} f\bigl(X^\varepsilon_s\bigr)\,\vd s \biggr] =\int
_{\RR} g_\alpha^\varepsilon(x,y)f(y)\,\vd y.
\]

The kernel $g_\varepsilon^\alpha$ satisfies the conditions
\begin{eqnarray*}
g^\varepsilon_\alpha(x,\varepsilon+)&=&g^\varepsilon_\alpha(x,
\varepsilon-),\qquad g^\varepsilon_\alpha(x,\varepsilon-)=g^\varepsilon_\alpha(x,
\varepsilon+),
\\
\nabla_y g^\varepsilon_\alpha(x,-\varepsilon-)&=&\kappa
\varepsilon\nabla_y g^\varepsilon_\alpha(x,-\varepsilon+),\\
\kappa\varepsilon\nabla_y g^\varepsilon_\alpha(x,
\varepsilon-)&=&\nabla_y g^\varepsilon_\alpha(x,
\varepsilon+),
\\
\nabla_yg^\varepsilon_\alpha(x,x+)-\nabla_y
g^\varepsilon_\alpha(x,x-)&=&2.
\end{eqnarray*}
With $\mu=\sqrt{2\alpha}$,
the coefficients
$A_\varepsilon$, $B_\varepsilon$, $C_\varepsilon$, $H_\varepsilon$ and
$F_\varepsilon$
are then expressed with the help of
\[
G_\varepsilon:= \bigl( 2{e^{4{{{ \mu}
\varepsilon}/{\sqrt{\kappa\varepsilon}}}}}\sqrt{\kappa\varepsilon
}+{e^{4{{{ \mu}\varepsilon}/{\sqrt{\kappa\varepsilon}}}}}
\kappa\varepsilon+2\sqrt{\kappa\varepsilon}-\kappa\varepsilon+{e^{4{
{{ \mu}\varepsilon}/{\sqrt{\kappa\varepsilon}}}}}-1
\bigr) {{ \mu}}.
\]
Since $\varepsilon\to0$,
$G_\varepsilon=4\sqrt{\kappa\varepsilon} (1+\frac{\mu}{\kappa}+\grandO
(\kappa\varepsilon) )\mu$.
After tedious computations,
\begin{eqnarray*}
A_\varepsilon(x)&=&- \bigl( {e^{4{{{ \mu}\varepsilon}/{\sqrt{\kappa
\varepsilon}}}}}\kappa\varepsilon-\kappa
\varepsilon-{e^{4{{{ \mu}\varepsilon}/{
\sqrt{\kappa\varepsilon}}}}}+1 \bigr) {e^{{ \mu} ( 2\varepsilon-x )
}}/G_\varepsilon
\mathop{\longrightarrow}_{\varepsilon\to0}A_0(x)\\
&:=& -\frac{e^{-\mu x}}{\kappa+\mu},
\\
B_\varepsilon(x)&=&B_0(x):=-\frac{e^{-{ \mu}x}}{\mu},
\\
C_\varepsilon(x)&=& -2\sinh(2\varepsilon/\sqrt{\kappa\varepsilon})\bigl(de^{-\mu x+2\mu\varepsilon} -e^{-\mu x+2\mu\varepsilon}
+de^{\mu x} +e^{-\mu x}
\bigr)
\\
&&{}\times {e^{{{2{\mu}\varepsilon}/{\sqrt{\kappa\varepsilon}}}}}
/G_\varepsilon
+4\sqrt{\kappa\varepsilon}e^{\mu x}\cosh(2\varepsilon/\sqrt{\kappa
\varepsilon}) {e^{{{2{\mu}\varepsilon}/{\sqrt{\kappa\varepsilon
}}}}} /G_\varepsilon\\
&&\hspace*{-12pt}\mathop{\longrightarrow}_{\varepsilon\to0}
C_0(x):=\frac{\kappa e^{\mu x}}{\mu(\kappa+\mu)},
\\
H_\varepsilon(x)&=&-2{e^{{{{ \mu} ( 3\varepsilon+\varepsilon\sqrt
{\kappa\varepsilon}-x
\sqrt{\kappa\varepsilon} ) }/{\sqrt{\kappa\varepsilon}}}}} ( 1+\sqrt
{\kappa\varepsilon} ) \sqrt
{\kappa\varepsilon} /G_\varepsilon\mathop{\longrightarrow}_{\varepsilon\to0}
H_0(x)\\
&:=&-\frac{\kappa e^{-\mu x}}{2\mu(\kappa+\mu)},
\\
E_\varepsilon(x)&=&-2{e^{{{{ \mu} ( \varepsilon+\varepsilon
\sqrt{\kappa\varepsilon}-x\sqrt{\kappa\varepsilon} ) }/{\sqrt{\kappa
\varepsilon}}}}} ( 1-\sqrt{\kappa\varepsilon} ) \sqrt{
\kappa\varepsilon} /G_\varepsilon\mathop{\longrightarrow}_{\varepsilon\to0}
H_0(x),
\\
F_\varepsilon(x)&=&-4\sqrt{\kappa\varepsilon} {e^{{
{{ \mu} ( 2\varepsilon+2\varepsilon\sqrt{\kappa\varepsilon}-x\sqrt
{\kappa\varepsilon} ) }/{\sqrt{\kappa\varepsilon}}}}}
/G_\varepsilon\mathop{\longrightarrow}_{\varepsilon\to0} F_0(x):=
\frac{-\kappa e^{-\mu x}}{\mu(\kappa+\mu)}\\
&=&-C_0(-x).
\end{eqnarray*}
Let $g_\alpha$ be the function
\[
g_\alpha(x,y):=\cases{ C_0(x)e^{-\mu y},&\quad$\mbox{if
}y>x$,\vspace*{2pt}
\cr
A_0(x)e^{-\mu y}+B_0(x)e^{\mu y},&\quad $
\mbox{if }y\in[0,x]$,\vspace*{2pt}
\cr
F_0(x)e^{\mu y},&\quad$\mbox{if }y<0$.}
\]
A similar work may be performed for $x<0$.
Thus, we easily obtain that $g^\varepsilon_\alpha(x,y)\longrightarrow
_{\varepsilon\to0}g_\alpha(x,y)$ converges to $g_\alpha$ and that
$g_\alpha$ is the density resolvent of the SNOB
by checking it satisfies the appropriate conditions at the interface.
\end{pf*}

%
\begin{proposition}
\label{prop-6}
Let $\fh^\varepsilon_0$ be the first hitting time of $0$
for $X^\varepsilon$.

Under $\PP_x$, $\fh^\varepsilon_0$ converges in
distribution to a random variable $\fk$
distributed as the lifetime of the EBM of parameter $\kappa$.
\end{proposition}

\begin{pf}
As in \cite{lejay-martinez,etore1},
we introduce $\Phi^\varepsilon(x)$ as the piecewise linear function
defined by
\[
\frac{\mathrm{d}\Phi^\varepsilon}{\mathrm{d} x}(x)=\cases{ 1/\sqrt
{\kappa\varepsilon},&\quad$\mbox{if }x\in[-
\varepsilon,\varepsilon]$,\vspace*{2pt}
\cr
1,&\quad $\mbox{otherwise}$.}
\]
Set $Y^\varepsilon=\Phi^\varepsilon(X^\varepsilon)$ so that
$Y^\varepsilon$
is solution to the SDE \cite{etore1,lejay-martinez}
\begin{eqnarray}
Y^\varepsilon_t=\Phi^\varepsilon(x)+B_t+
\theta^\varepsilon L^{y_\varepsilon}_t\bigl(Y^\varepsilon\bigr) -
\theta^\varepsilon L^{-y_\varepsilon}_t\bigl(Y^\varepsilon\bigr)
\nonumber\\
\eqntext{\mbox{with } \displaystyle\theta^\varepsilon=\frac{1-\sqrt{\kappa\varepsilon
}}{1+\sqrt{\kappa\varepsilon}} \mbox{ and }y^\varepsilon:=\Phi^\varepsilon(\varepsilon)=\sqrt{\frac{\varepsilon
}{\kappa}}.}
\end{eqnarray}
The infinitesimal generator of $Y^\varepsilon$
is $\cA^\varepsilon:=\frac{1}{2}\triangle$ whose domain contains as a dense
subset [it is similar to the discussion on $D^\varepsilon$ in \eqref{eq:Deps}]
\[
\lleft\{ f\in\cC^2\bigl((-\infty,-y_\varepsilon)\cup
(-y_\varepsilon,y_\varepsilon)\cup(y_\varepsilon,\infty)\bigr) \middle|
\\
\begin{array} {l} f,f''\in\uL^2(\RR),
\\
{}[f]\bigl(\pm y^\varepsilon\bigr)=0,
\\
\bigl(1-\theta^\varepsilon\bigr)f'\bigl(y^\varepsilon-
\bigr)\\
\qquad=\bigl(1+\theta^\varepsilon\bigr)f'\bigl(y^\varepsilon+
\bigr),
\\
\bigl(1+\theta^\varepsilon\bigr)f'\bigl(-y^\varepsilon-
\bigr)\\
\qquad=\bigl(1-\theta^\varepsilon\bigr)f'\bigl(-y^\varepsilon+
\bigr),
\end{array}
\rright\}.
\]

From now, we assume for the sake of simplicity that $x>0$.

The hitting time $\fh^\varepsilon_0$ is also
the first hitting time of zero by $Y^\varepsilon$.
Since by symmetry $\psi(-x,\alpha)=\psi(x,\alpha)$ for any $x\geq0$,
we consider
only that $x\geq0$.

Since the Feynman--Kac formula is valid for
the process $Y^\varepsilon$, $\psi^\varepsilon(x,\alpha):=\EE
_x[e^{-\alpha\fh^\varepsilon_0}]$
is solution to
\[
\cases{ \frac{1}{2}\triangle\psi^\varepsilon(x,\alpha)=\alpha
\psi^\varepsilon(x,\alpha),&\quad$\mbox{for }x\neq y^\varepsilon,$\vspace*{2pt}
\cr
\psi^\varepsilon(0,\alpha)=1,&\vspace*{2pt}
\cr
\psi^\varepsilon
\bigl(y^\varepsilon-,\alpha\bigr)=\psi^\varepsilon\bigl(y^\varepsilon+,
\alpha\bigr),&\vspace*{2pt}
\cr
\bigl(1-\theta^\varepsilon\bigr)
\nabla_x\psi^\varepsilon\bigl(y^\varepsilon-,\alpha\bigr) =
\bigl(1+\theta^\varepsilon\bigr)\nabla_x\psi^\varepsilon
\bigl(y^\varepsilon+,\alpha\bigr). }
\]
Hence, $\psi^\varepsilon(x,\alpha)$ is sought as
\[
\psi^\varepsilon(x,\alpha)=\cases{ \gamma^\varepsilon\exp(-\sqrt{2
\alpha}x),&\quad$\mbox{if }x>y^\varepsilon$,\vspace*{2pt}
\cr
\cos(\sqrt{2\alpha}x)+
\beta^\varepsilon\sin(\sqrt{2\alpha}x),&\quad$\mbox{if }x\in\bigl
[0,y^\varepsilon
\bigr]$. }
\]
After some computations,
\[
\beta^\varepsilon=\frac{-\cos(\sqrt{2\alpha} y^\varepsilon)+\sqrt{\kappa
\varepsilon}\sin(\sqrt{2\alpha} y^\varepsilon)} {
\sin(\sqrt{2\alpha} y^\varepsilon)+\sqrt{\kappa\varepsilon}\cos(\sqrt
{2\alpha} y^\varepsilon)} \quad\mbox{and}\quad\sqrt{\varepsilon}
\beta^\varepsilon\mathop{\sim}_{\varepsilon\to0}\frac{-\sqrt{\kappa}}{\kappa
+\sqrt{2\alpha}}.
\]
Besides,
\[
\gamma^\varepsilon=e^{\sqrt{2\alpha}y^\varepsilon}\sqrt{\kappa
\varepsilon} \bigl(
\beta^\varepsilon\cos\bigl(\sqrt{2\alpha}y^\varepsilon\bigr)-
\beta^\varepsilon\sin\bigl(\sqrt{2\alpha}y^\varepsilon\bigr) \bigr)
\mathop{\sim}_{\varepsilon\to0}\frac{\kappa}{\kappa+\sqrt{2\alpha}}.
\]
Hence, for any $x>0$,
%
%
\begin{equation}
\label{eq-time1} \psi^\varepsilon(x,\alpha)\mathop{\longrightarrow}_{\varepsilon\to
0}\psi(x,
\alpha):=\frac{\kappa}{\kappa+\sqrt{2\alpha}} e^{-\sqrt{2\alpha}x}
\end{equation}
with $\psi$ defined by \eqref{eq:ebm:2}.

This proves that under $\PP_x$, $\fh^\varepsilon_0$ converges to
a random variable $\fk$ whose Laplace transform
is $\psi(x,\alpha)$ under $\PP_x$. This random variable $\fk$
is then the lifetime of an EBM.
\end{pf}

\begin{pf*}{Proof of Theorem~\ref{thm:main}}
Using the properties of the resolvent,
for $\alpha>0$ and a bounded, measurable function $f$,
\begin{eqnarray*}
G^\varepsilon_\alpha f(x) &:=&\EE_x \biggl[\int
_0^{+\infty} e^{-\alpha t}f\bigl(X^\varepsilon_s
\bigr)\,\vd s \biggr]\\
& =&R_\alpha^\varepsilon f(x) +\EE_x
\bigl[e^{-\alpha\fh_0^\varepsilon}\bigr] \frac{1}{2} \bigl(
G^\varepsilon_\alpha
f(\varepsilon) +G^\varepsilon_\alpha f(-\varepsilon) \bigr)
\end{eqnarray*}
with
\[
R_\alpha^\varepsilon f(x):= \EE_x \biggl[\int
_0^{\fh_0^\varepsilon} e^{-\alpha t}f\bigl(X^\varepsilon_s
\bigr)\,\vd s \biggr].
\]
Since $\psi^\varepsilon(x,\alpha)=\psi^\varepsilon(-x,\alpha)$,
\[
G_\alpha^\varepsilon f(x)=R_\varepsilon^\alpha f(x) +
\frac{\psi^\varepsilon(x,\alpha)}{1-\psi^\varepsilon(\varepsilon,\alpha
)}\frac{R_\varepsilon^\alpha f(\varepsilon)+R_\varepsilon^\alpha
f(-\varepsilon)}{2}.
\]

For the sake of simplicity, we assume that $x>0$.
Using the symmetry properties of $\cL^\varepsilon$,
\[
R_\alpha^\varepsilon f(x) =\int_0^{+\infty}
\bigl(g_\alpha^\varepsilon(x,y)-g^\varepsilon_\alpha(x,-y)
\bigr)f(y)\,\vd y.
\]
But
\[
g_\alpha^\varepsilon(x,y)-g^\varepsilon_\alpha(x,-y)
\mathop{\longrightarrow}_{\varepsilon\to0} g_\alpha(x,y)-g_\alpha(x,-y)=
\gel_\alpha(x,y),
\]
where $\gel_\alpha(x,y)$ is the resolvent density of the EBM.
Thus, $R_\varepsilon^\alpha f(x)\longrightarrow_{\varepsilon\to0} \Gel
_\alpha f(x)$
for any $x>0$.
It is also easily obtained that
\[
R_\varepsilon^\alpha f(\varepsilon)\mathop{\longrightarrow}_{\varepsilon\to0}
\Gel_\alpha f(0+)\quad \mbox{and} \quad R_\varepsilon^\alpha f(-
\varepsilon)\mathop{\longrightarrow}_{\varepsilon\to0} \Gel_\alpha f(0-).
\]
Using \eqref{eq:1bis} and \eqref{eq-time1},
$G_\alpha^\varepsilon f(x)\longrightarrow_{\varepsilon\to0}
G_\alpha f(x)$. The Trotter--Kato theorem
(see, e.g., \cite{kato}, Theorem IX.2.16, page 504)
and the Markov property imply the convergence
in finite-dimensional distributions
of $Z^\varepsilon$ to $X$ under $\PP_x$ for $x\geq0$.
By symmetry, this could be extended to $x\leq0$.

The only remaining point of the tightness.
When away from $[-\varepsilon,\varepsilon]$,
$X^\varepsilon$ behaves like a Brownian motion.
Hence, for $0\leq s\leq t\leq T$, let us set
$\ff(s,t):=\inf\{u>s;|X^\varepsilon_u|=\varepsilon\}$
with possibly $\ff(s,t)=+\infty$
and $\fl(s,t):=\sup\{u<t;|X^\varepsilon_u|=\varepsilon\}$
with possibly $\fl(s,t)=-\infty$.

If $\ff(s,t)\geq t$ and $\fl(s,t)\leq s$, then for $\delta<1/2$,
there exists an integrable random variable $C(\omega)$
such that
$|X^\varepsilon_t(\omega)-X^\varepsilon_s(\omega)|\leq C(\omega
)(t-s)^\delta$
for any $0\leq s\leq t\leq T$.

If $\ff(s,t)\leq t$ and $\fl(s,t)\leq s$, then
\[
\bigl|X^\varepsilon_t-X^\varepsilon_s\bigr|\leq
\bigl|X^\varepsilon_{\ff(s,t)}-X^\varepsilon_s\bigr|
+\bigl|X^\varepsilon_t-X^\varepsilon_{\ff(s,t)}\bigr| \leq
C(t-s)^\beta+2\varepsilon
\]
since $X^\varepsilon_t$ belongs to $[-\varepsilon,\varepsilon]$.
A similar analysis could be carried for the other cases,
which means that for some integrable random variable $C$,
\[
\sup_{|t-s|<\delta}\bigl|X^\varepsilon_t-X^\varepsilon_s\bigr|
\leq C\delta^\beta+2\varepsilon.
\]
This proves that $(Z^\varepsilon)_{\varepsilon>0}$ is tight is
the space $\mathcal{D}([0,T];\RR)$ of discontinuous functions with the
Skorohod topology
(see, e.g., \cite{billingsley})
and then on $\mathcal{D}([0,T];\GG)$.
Hence, we easily deduce the convergence of $Z^\varepsilon$
to the SNOB in $\mathcal{D}([0,T];\GG)$.
\end{pf*}

\section{Simulation of the SNOB}
\label{sec:simulation}

It is easy to simulate a discretized process $X$ in the same way
it is easy to simulate the Brownian motion.
Following Proposition~\ref{prop:representation},
we draw a random variate
with density $p(\delta t,x,\cdot)$ when $x$ is close enough to~$0$.

For this, we use a Brownian bridge technique to check
if the process reaches $0\pm$ before $\delta t$ (see, e.g., \cite{baldi}
and \cite{lejay-pichot}, Section B.2, for an example of application
and further references). This involve the inverse Gaussian distribution
$\cIG(\lambda,\mu)$ whose density is
$r_{\mu,\lambda}(x)=\sqrt{\frac{\lambda}{2\pi x^3}}
\exp(\frac{-\lambda(x-\mu)^2}{2\mu^2 x} )$.
Random variates with $\cIG$ distribution could be simulated by the
methods proposed in
\cite{devroye86a},  page 148 and~\cite{michael}.

We simulate the local time using the following representation
under $\PP_0$ \cite{lepingle93a,lepingle95a}:
\[
\bigl(L_t^0(B),|B_t|\bigr)\equaldist(\fl,
\fl-H) \qquad\mbox{where }\fl:=\tfrac{1}{2}\bigl(H+\sqrt{V+H^2}\bigr)
\]
with $H\sim\cN(0,t)$ and $V\sim\exp(1/2t)$ independent from $H$.

The generic algorithm to simulate
the process at time $\delta t$ when at point $x$ at time~$0$ is the following:
\begin{enumerate}[1.]
\item[1.] Set $y:=x+\sqrt{\delta t}G$ with $G$ a random variate whose distribution
is $\cN(0,1)$.
\item[2.] If $|x|\geq4\sqrt{\delta t}$, then return $y$ (here, we neglect
the exponentially small probability that the process crosses
$0$ between the times $0$ and $\delta t$).
\item[3.] If $xy>0$, then decide with probability $\exp(-2|xy|/\delta t)$ if
the path $X$ has crossed $0$.
\begin{itemize}
\item If no crossing occurs, then return $y$.
\item If a crossing occurs, draw $\fg\sim\cIG
(|x|/|y|,x^2/2\delta t)$, so that $\fz:=\delta t \fg/(1+\fg)$ is a
realization of the first hitting time of $0$
for a Brownian bridge with $B_0=x$ and $B_{\delta t}=y$.
Then go the step 5.
\end{itemize}
\item[4.] If $xy<0$, then draw $\fg\sim\cIG(-|x|/|y|,x^2/2\delta t)$
and set $\fz:=\delta t \fg/(1+\fg)$, the first time the Brownian bridge
reaches $0$.
Go to step 5.
\item[5.]
Set $\fr:=\delta t-\fz$.
For two independent random variates
$H\sim\cN(0,\fr)$ and $V\sim\exp(1/2\fr)$,
set $\fl:=(H+\sqrt{V+H^2})/2$.
\item[6.]
For $U\sim\cU(0,1)$ independent from $V$ and $H$,
set $\fs:=\sgn(x)$ if $\exp(-\kappa\fl)\geq2U-1$.
Otherwise, set $\fs:=-\sgn(x)$.
\item[7.] Return $\fs(\fl-H)$.
\end{enumerate}

An application to the estimation of a macroscopic estimation parameter
in the context of a simplified problem related to brain imaging may be
found in \cite{lejay-2015}.
The results are satisfactory, unless $\kappa$ is too small due to a problem
of rare event simulation.

\section*{Acknowledgements} The author is indebted to
Jing-Rebecca Li and Denis Grebenkov for having proposed
this research and interesting discussions about it.
The author also wishes to thank warmly his wife, Claire Nivlet, for
having suggested the name of the process.



\begin{thebibliography}{36}

\bibitem{andrews}
%
\begin{barticle}[pbm]
\bauthor{\bsnm{Andrews},~\bfnm{Steven~S.}\binits{S.~S.}}
(\byear{2009}).
\btitle{Accurate particle-based simulation of adsorption, desorption
and partial transmission}.
\bjournal{Phys. Biol.}
\bvolume{6}
\bpages{046015}.
\bid{doi={10.1088/1478-3975/6/4/046015}, issn={1478-3975},
mid={NIHMS170155}, pii={S1478-3975(09)23496-0}, pmcid={2847898},
pmid={19910670}}
\end{barticle}
%
\bptok{imsref}%
\endbibitem

\bibitem{baldi}
%
\begin{barticle}[mr]
\bauthor{\bsnm{Baldi},~\bfnm{Paolo}\binits{P.}}
(\byear{1995}).
\btitle{Exact asymptotics for the probability of exit from a domain and
applications to simulation}.
\bjournal{Ann. Probab.}
\bvolume{23}
\bpages{1644--1670}.
\bid{issn={0091-1798}, mr={1379162}}
\end{barticle}
%
\bptok{imsref}%
\endbibitem

\bibitem{bass08a}
%
\begin{barticle}[mr]
\bauthor{\bsnm{Bass},~\bfnm{Richard~F.}\binits{R.~F.}},
\bauthor{\bsnm{Burdzy},~\bfnm{Krzysztof}\binits{K.}} \AND
\bauthor{\bsnm{Chen},~\bfnm{Zhen-Qing}\binits{Z.-Q.}}
(\byear{2008}).
\btitle{On the {R}obin problem in fractal domains}.
\bjournal{Proc. Lond. Math. Soc. (3)}
\bvolume{96}
\bpages{273--311}.
\bid{doi={10.1112/plms/pdm045}, issn={0024-6115}, mr={2396121}}
\end{barticle}
%
\bptok{imsref}%
\endbibitem

\bibitem{bass}
%
\begin{barticle}[mr]
\bauthor{\bsnm{Bass},~\bfnm{Richard~F.}\binits{R.~F.}} \AND
\bauthor{\bsnm{Chen},~\bfnm{Zhen-Qing}\binits{Z.-Q.}}
(\byear{2003}).
\btitle{Brownian motion with singular drift}.
\bjournal{Ann. Probab.}
\bvolume{31}
\bpages{791--817}.
\bid{doi={10.1214/aop/1048516536}, issn={0091-1798}, mr={1964949}}
\end{barticle}
%
\bptok{imsref}%
\endbibitem

\bibitem{billingsley}
%
\begin{bbook}[mr]
\bauthor{\bsnm{Billingsley},~\bfnm{Patrick}\binits{P.}}
(\byear{1999}).
\btitle{Convergence of Probability Measures},
\bedition{2nd} ed.
\bpublisher{Wiley},
\blocation{New York}.
\bid{doi={10.1002/9780470316962}, mr={1700749}}
\end{bbook}
%
\bptok{imsref}%
\endbibitem

\bibitem{bobrowski}
%
\begin{barticle}[mr]
\bauthor{\bsnm{Bobrowski},~\bfnm{Adam}\binits{A.}}
(\byear{2012}).
\btitle{From diffusions on graphs to {M}arkov chains via asymptotic
state lumping}.
\bjournal{Ann. Henri Poincar\'e}
\bvolume{13}
\bpages{1501--1510}.
\bid{doi={10.1007/s00023-012-0158-z}, issn={1424-0637}, mr={2966471}}
\end{barticle}
%
\bptok{imsref}%
\endbibitem

\bibitem{devroye86a}
%
\begin{bbook}[mr]
\bauthor{\bsnm{Devroye},~\bfnm{Luc}\binits{L.}}
(\byear{1986}).
\btitle{Nonuniform Random Variate Generation}.
\bpublisher{Springer},
\blocation{New York}.
\bid{doi={10.1007/978-1-4613-8643-8}, mr={0836973}}
\end{bbook}
%
\bptok{imsref}%
\endbibitem

\bibitem{erban}
%
\begin{barticle}[author]
\bauthor{\bsnm{Erban},~\bfnm{R.}\binits{R.}} \AND
\bauthor{\bsnm{Chapman},~\bfnm{S.~J.}\binits{S.~J.}}
(\byear{2007}).
\btitle{Reactive boundary conditons for stochastic simulation of
reaction--diffusion processes}.
\bjournal{Phys. Biol.}
\bvolume{4}
\bpages{16--28}.
\end{barticle}
%
\bptok{imsref}%
\endbibitem

\bibitem{etore1}
%
\begin{barticle}[mr]
\bauthor{\bsnm{{\'E}tor{\'e}},~\bfnm{Pierre}\binits{P.}}
(\byear{2006}).
\btitle{On random walk simulation of one-dimensional diffusion
processes with discontinuous coefficients}.
\bjournal{Electron. J. Probab.}
\bvolume{11}
\bpages{249--275 (electronic)}.
\bid{doi={10.1214/EJP.v11-311}, issn={1083-6489}, mr={2217816}}
\end{barticle}
%
\bptok{imsref}%
\endbibitem

\bibitem{feller}
%
\begin{barticle}[mr]
\bauthor{\bsnm{Feller},~\bfnm{William}\binits{W.}}
(\byear{1954}).
\btitle{Diffusion processes in one dimension}.
\bjournal{Trans. Amer. Math. Soc.}
\bvolume{77}
\bpages{1--31}.
\bid{issn={0002-9947}, mr={0063607}}
\end{barticle}
%
\bptok{imsref}%
\endbibitem

\bibitem{fiermans}
%
\begin{barticle}[pbm]
\bauthor{\bsnm{Fieremans},~\bfnm{Els}\binits{E.}},
\bauthor{\bsnm{Novikov},~\bfnm{Dmitry~S.}\binits{D.~S.}},
\bauthor{\bsnm{Jensen},~\bfnm{Jens~H.}\binits{J.~H.}} \AND
\bauthor{\bsnm{Helpern},~\bfnm{Joseph~A.}\binits{J.~A.}}
(\byear{2010}).
\btitle{Monte Carlo study of a two-compartment exchange model of diffusion}.
\bjournal{NMR Biomed.}
\bvolume{23}
\bpages{711--724}.
\bid{doi={10.1002/nbm.1577}, issn={1099-1492}, mid={NIHMS248852},
pmcid={2997614}, pmid={20882537}}
\end{barticle}
%
\bptok{imsref}%
\endbibitem

\bibitem{freidlin}
%
\begin{barticle}[mr]
\bauthor{\bsnm{Freidlin},~\bfnm{Mark~I.}\binits{M.~I.}} \AND
\bauthor{\bsnm{Wentzell},~\bfnm{Alexander~D.}\binits{A.~D.}}
(\byear{1993}).
\btitle{Diffusion processes on graphs and the averaging principle}.
\bjournal{Ann. Probab.}
\bvolume{21}
\bpages{2215--2245}.
\bid{issn={0091-1798}, mr={1245308}}
\end{barticle}
%
\bptok{imsref}%
\endbibitem

\bibitem{freidlin2}
%
\begin{barticle}[mr]
\bauthor{\bsnm{Freidlin},~\bfnm{Mark~I.}\binits{M.~I.}} \AND
\bauthor{\bsnm{Wentzell},~\bfnm{Alexander~D.}\binits{A.~D.}}
(\byear{1994}).
\btitle{Random perturbations of {H}amiltonian systems}.
\bjournal{Mem. Amer. Math. Soc.}
\bvolume{109}
\bpages{viii+82}.
\bid{doi={10.1090/memo/0523}, issn={0065-9266}, mr={1201269}}
\bptnote{check pages}%
\end{barticle}
%
\bptok{imsref}%
\endbibitem

\bibitem{gallavotti}
\begin{barticle}[mr]
\bauthor{\bsnm{Gallavotti},~\bfnm{G.}\binits{G.}} \AND
\bauthor{\bsnm{McKean},~\bfnm{H.~P.}\binits{H.~P.}}
(\byear{1972}).
\btitle{Boundary conditions for the heat equation in a
several-dimensional region}.
\bjournal{Nagoya Math. J.}
\bvolume{47}
\bpages{1--14}.
\bid{issn={0027-7630}, mr={0317658}}
\end{barticle}
\bptok{imsref}\endbibitem\

\bibitem{grebenkov}
%
\begin{bincollection}[mr]
\bauthor{\bsnm{Grebenkov},~\bfnm{Denis~S.}\binits{D.~S.}}
(\byear{2006}).
\btitle{Partially reflected {B}rownian motion: A stochastic approach to
transport phenomena}.
In \bbooktitle{Focus on Probability Theory}
(\beditor{\bfnm{L.~R.}\binits{L.~R.}~\bsnm{Velle}}, ed.)
\bpages{135--169}.
\bpublisher{Nova Sci. Publ.},
\blocation{New York}.
\bid{mr={2553673}}
\end{bincollection}
%
\bptok{imsref}%
\endbibitem

\bibitem{griego}
%
\begin{barticle}[mr]
\bauthor{\bsnm{Griego},~\bfnm{Richard~J.}\binits{R.~J.}} \AND
\bauthor{\bsnm{Moncayo},~\bfnm{Alberto}\binits{A.}}
(\byear{1970}).
\btitle{Random evolutions and piecing out of {M}arkov processes}.
\bjournal{Bol. Soc. Mat. Mexicana (2)}
\bvolume{15}
\bpages{22--29}.
\bid{mr={0365723}}
\end{barticle}
%
\bptok{imsref}%
\endbibitem

\bibitem{ikeda}
%
\begin{barticle}[mr]
\bauthor{\bsnm{Ikeda},~\bfnm{Nobuyuki}\binits{N.}},
\bauthor{\bsnm{Nagasawa},~\bfnm{Masao}\binits{M.}} \AND
\bauthor{\bsnm{Watanabe},~\bfnm{Shinzo}\binits{S.}}
(\byear{1966}).
\btitle{A construction of {M}arkov processes by piecing out}.
\bjournal{Proc. Japan Acad.}
\bvolume{42}
\bpages{370--375}.
\bid{issn={0021-4280}, mr={0202197}}
\end{barticle}
%
\bptok{imsref}%
\endbibitem

\bibitem{ito-mckean}
%
\begin{bbook}[mr]
\bauthor{\bsnm{It{\^o}},~\bfnm{Kiyoshi}\binits{K.}} \AND
\bauthor{\bsnm{McKean},~\bfnm{Henry~P.}\binits{H.~P.} \bsuffix{Jr.}}
(\byear{1996}).
\btitle{Diffusion Processes and Their Sample Paths}.
\bpublisher{Springer},
\blocation{Berlin}.
\end{bbook}
%
\bptok{imsref}%
\endbibitem

\bibitem{karlin}
%
\begin{bbook}[mr]
\bauthor{\bsnm{Karlin},~\bfnm{Samuel}\binits{S.}} \AND
\bauthor{\bsnm{Taylor},~\bfnm{Howard~M.}\binits{H.~M.}}
(\byear{1981}).
\btitle{A Second Course in Stochastic Processes}.
\bpublisher{Academic Press},
\blocation{New York}.
\bid{mr={0611513}}
\end{bbook}
%
\bptok{imsref}%
\endbibitem

\bibitem{kato}
%
\begin{bbook}[mr]
\bauthor{\bsnm{Kato},~\bfnm{Tosio}\binits{T.}}
(\byear{1995}).
\btitle{Perturbation Theory for Linear Operators}.
\bpublisher{Springer},
\blocation{Berlin}.
\bid{mr={1335452}}
\end{bbook}
%
\bptok{imsref}%
\endbibitem

\bibitem{kopytko}
%
\begin{barticle}[mr]
\bauthor{\bsnm{Kopytko},~\bfnm{Bohdan~I.}\binits{B.~I.}} \AND
\bauthor{\bsnm{Portenko},~\bfnm{Mykola~I.}\binits{M.~I.}}
(\byear{2009}).
\btitle{The problem of pasting together two diffusion processes and
classical potentials}.
\bjournal{Theory Stoch. Process.}
\bvolume{15}
\bpages{126--139}.
\bid{issn={0321-3900}, mr={2598532}}
\end{barticle}
%
\bptok{imsref}%
\endbibitem

\bibitem{lejay-sbm}
%
\begin{barticle}[mr]
\bauthor{\bsnm{Lejay},~\bfnm{Antoine}\binits{A.}}
(\byear{2006}).
\btitle{On the constructions of the skew {B}rownian motion}.
\bjournal{Probab. Surv.}
\bvolume{3}
\bpages{413--466}.
\bid{doi={10.1214/154957807000000013}, issn={1549-5787}, mr={2280299}}
\end{barticle}
%
\bptok{imsref}%
\endbibitem

\bibitem{lejay-2015}
%
\begin{bmisc}[author]
\bauthor{\bsnm{Lejay},~\bfnm{A.}\binits{A.}}
(\byear{2015}).
\bhowpublished{Estimation of the mean residence time in cells
surrounded by semi-permeable membranes by a Monte Carlo method.
Research report No. Inria, RR-8709}.
\end{bmisc}
%
\bptok{imsref}%
\endbibitem

\bibitem{lejay-martinez}
%
\begin{barticle}[mr]
\bauthor{\bsnm{Lejay},~\bfnm{Antoine}\binits{A.}} \AND
\bauthor{\bsnm{Martinez},~\bfnm{Miguel}\binits{M.}}
(\byear{2006}).
\btitle{A scheme for simulating one-dimensional diffusion processes
with discontinuous coefficients}.
\bjournal{Ann. Appl. Probab.}
\bvolume{16}
\bpages{107--139}.
\bid{doi={10.1214/105051605000000656}, issn={1050-5164}, mr={2209338}}
\end{barticle}
%
\bptok{imsref}%
\endbibitem

\bibitem{lejay-pichot}
%
\begin{barticle}[mr]
\bauthor{\bsnm{Lejay},~\bfnm{Antoine}\binits{A.}} \AND
\bauthor{\bsnm{Pichot},~\bfnm{G{\'e}raldine}\binits{G.}}
(\byear{2012}).
\btitle{Simulating diffusion processes in discontinuous media: A
numerical scheme with constant time steps}.
\bjournal{J. Comput. Phys.}
\bvolume{231}
\bpages{7299--7314}.
\bid{doi={10.1016/j.jcp.2012.07.011}, issn={0021-9991}, mr={2969713}}
\end{barticle}
%
\bptok{imsref}%
\endbibitem

\bibitem{lepingle93a}
%
\begin{barticle}[mr]
\bauthor{\bsnm{L{\'e}pingle},~\bfnm{Dominique}\binits{D.}}
(\byear{1993}).
\btitle{Un sch\'ema d'{E}uler pour \'equations diff\'erentielles
stochastiques r\'efl\'echies}.
\bjournal{C. R. Acad. Sci. Paris S\'er. I Math.}
\bvolume{316}
\bpages{601--605}.
\bid{issn={0764-4442}, mr={1212213}}
\end{barticle}
%
\bptok{imsref}%
\endbibitem

\bibitem{lepingle95a}
%
\begin{barticle}[mr]
\bauthor{\bsnm{L{\'e}pingle},~\bfnm{D.}\binits{D.}}
(\byear{1995}).
\btitle{Euler scheme for reflected stochastic differential equations}.
\bjournal{Math. Comput. Simulation}
\bvolume{38}
\bpages{119--126}.
\bid{doi={10.1016/0378-4754(93)E0074-F}, issn={0378-4754}, mr={1341164}}
\end{barticle}
%
\bptok{imsref}%
\endbibitem

\bibitem{legall}
%
\begin{bincollection}[mr]
\bauthor{\bsnm{Le Gall},~\bfnm{J.-F.}\binits{J.-F.}}
(\byear{1984}).
\btitle{One-dimensional stochastic differential equations involving the
local times of the unknown process}.
In \bbooktitle{Stochastic Analysis and Applications ({S}wansea, 1983)}.
\bseries{Lecture Notes in Math.}
\bvolume{1095}
\bpages{51--82}.
\bpublisher{Springer},
\blocation{Berlin}.
\bid{doi={10.1007/BFb0099122}, mr={0777514}}
\end{bincollection}
%
\bptok{imsref}%
\endbibitem

\bibitem{meyer}
%
\begin{barticle}[mr]
\bauthor{\bsnm{Meyer},~\bfnm{P.~A.}\binits{P.~A.}}
(\byear{1975}).
\btitle{Renaissance, recollements, m\'elanges, ralentissement de
processus de {M}arkov}.
\bjournal{Ann. Inst. Fourier (Grenoble)}
\bvolume{25}
\bpages{465--497}.
\bid{issn={0373-0956}, mr={0415784}}
\end{barticle}
%
\bptok{imsref}%
\endbibitem

\bibitem{michael}
%
\begin{barticle}[author]
\bauthor{\bsnm{Michael},~\bfnm{J.~R.}\binits{J.~R.}},
\bauthor{\bsnm{Shucany},~\bfnm{W.~R.}\binits{W.~R.}} \AND
\bauthor{\bsnm{Haas},~\bfnm{R.~W.}\binits{R.~W.}}
(\byear{1976}).
\btitle{Generating random variates using transformations with multiple roots}.
\bjournal{Amer. Statist.}
\bvolume{30}
\bpages{88--90}.
\end{barticle}
%
\bptok{imsref}%
\endbibitem

\bibitem{papanicolaou}
%
\begin{barticle}[mr]
\bauthor{\bsnm{Papanicolaou},~\bfnm{Vassilis~G.}\binits{V.~G.}}
(\byear{1990}).
\btitle{The probabilistic solution of the third boundary value problem
for second order elliptic equations}.
\bjournal{Probab. Theory Related Fields}
\bvolume{87}
\bpages{27--77}.
\bid{doi={10.1007/BF01217746}, issn={0178-8051}, mr={1076956}}
\end{barticle}
%
\bptok{imsref}%
\endbibitem

\bibitem{portenko}
%
\begin{barticle}[mr]
\bauthor{\bsnm{Portenko},~\bfnm{N.~I.}\binits{N.~I.}}
(\byear{2000}).
\btitle{A probabilistic representation for the solution to one problem
of mathematical physics}.
\bjournal{Ukrainian Math. J.}
\bvolume{52}
\bpages{1457--1469}.
\bptnote{check pages}%
\end{barticle}
%
\bptok{imsref}%
\endbibitem

\bibitem{sanchez}
%
\begin{bbook}[author]
\bauthor{\bsnm{S{\'a}nchez-Palencia},~\bfnm{E.}\binits{E.}}
(\byear{1980}).
\btitle{Non-Homogeneous Media and Vibration Theory}.
\bseries{Lecture Notes in Physics}
\bvolume{127}.
\bpublisher{Springer},
\blocation{Berlin}.
\end{bbook}
%
\bptok{imsref}%
\endbibitem

\bibitem{siegrist}
%
\begin{barticle}[mr]
\bauthor{\bsnm{Siegrist},~\bfnm{Kyle}\binits{K.}}
(\byear{1981}).
\btitle{Random evolution processes with feedback}.
\bjournal{Trans. Amer. Math. Soc.}
\bvolume{265}
\bpages{375--392}.
\bid{doi={10.2307/1999740}, issn={0002-9947}, mr={0610955}}
\end{barticle}
%
\bptok{imsref}%
\endbibitem

\bibitem{singer}
\begin{barticle}[mr]
\bauthor{\bsnm{Singer},~\bfnm{A.}\binits{A.}},
\bauthor{\bsnm{Schuss},~\bfnm{Z.}\binits{Z.}},
\bauthor{\bsnm{Osipov},~\bfnm{A.}\binits{A.}} \AND
\bauthor{\bsnm{Holcman},~\bfnm{D.}\binits{D.}}
(\byear{2007/08}).
\btitle{Partially reflected diffusion}.
\bjournal{SIAM J. Appl. Math.}
\bvolume{68}
\bpages{844--868}.
\bid{doi={10.1137/060663258}, issn={0036-1399}, mr={2375298}}
\end{barticle}
\bptok{imsref}\endbibitem\

\bibitem{stroock}
%
\begin{bincollection}[mr]
\bauthor{\bsnm{Stroock},~\bfnm{Daniel~W.}\binits{D.~W.}}
(\byear{1988}).
\btitle{Diffusion semigroups corresponding to uniformly elliptic
divergence form operators}.
In \bbooktitle{S\'eminaire de {P}robabilit\'es, {XXII}}.
\bseries{Lecture Notes in Math.}
\bvolume{1321}
\bpages{316--347}.
\bpublisher{Springer},
\blocation{Berlin}.
\bid{doi={10.1007/BFb0084145}, mr={0960535}}
\end{bincollection}
%
\bptok{imsref}%
\endbibitem

\end{thebibliography}

%





\printaddresses
\end{document}